# Fundamental trigonometric interpolation and approximating polynomials and splines

Denysiuk V.P.
Dr. of Phys-Math. sciences, Professor, Kiev, Ukraine
National Aviation University
kvomden@nau.edu.ua

## Annotation

The paper deals with two fundamental types of trigonometric polynomials and splines on uniform grids, which allow us to construct interpolation approximations that depend linearly on the values of the interpolated function. Fundamental on the same grids are trigonometric LS polynomials and LS splines, which allow us to construct the approximation of functions by the least squares method and also depend linearly on the values of the approximate function. Trigonometric splines were considered as splines with polynomial analogues in some cases. The material is illustrated in graphs. The considered fundamental trigonometric polynomials and splines can be recommended for use in many problems related to approximation of functions, in particular, mathematical modeling of signals and processes.

## Keywords:

## Introduction

Consideration of systems of fundamental functions on a grid is advisable to begin with consideration of the interpolation problem.

Let for a segment $[0,T]$, some grid is specified $\Delta_N$, $\Delta_N = \{t_i\}_{i=0}^{N}$, $0 \le t_0 < t_1 < \ldots < t_N \le 1$. Let also on segment $[0,T]$ function is specified $f(t)$, and the exact values are known $f(t_j) = f_j$, $(j = 0,1,\ldots,N)$, this feature in the grid nodes $\Delta_N$. We need to construct a generalized polynomial

$$\Phi_N(t) = c_0\varphi_0(t) + c_1\varphi_1(t) + \ldots + c_N\varphi_N(t)$$

on a system of linearly independent functions $\varphi_0(t), \varphi_1(t), \ldots \varphi_N(t)$, which depends on $N+1$- st parameters $c_0, c_1, \ldots, c_N$ and satisfies the conditions

$$\Phi_N(x_j) = f_i, \quad j = 0,1,\ldots,N.$$

The formulated problem is called the interpolation problem [1,2], and a polynomial $\Phi_N(x)$ is called a generalized interpolation polynomial.

In many cases, in interpolation problems, the system of functions $\varphi_0(t), \varphi_1(t), \ldots, \varphi_N(t)$ on grid $\Delta_N$ it is advisable to choose such coefficients which $c_k$, $(k = 0,1,\ldots,N)$, generalized interpolation polynomial $\Phi_N(x)$ represented the value of the interpolated function $f(t)$ in the nodes of this grid. In this case, the polynomial $\Phi_N(t)$ takes shape

$$\Phi_N(t) = f_0\varphi_0(t) + f_1\varphi_1(t) + \ldots + f_N\varphi_N(t). \tag{1}$$

It is clear that in this case $\varphi(t)$ the functions must be correlated

$$\varphi_j(t_i) = \begin{cases} 1, & j = i; \\ 0, & j \neq i. \end{cases} \tag{2}$$

Functions for which the grid there are relations (2), called fundamental functions on this grid [3]. It is clear that the area of definition of fundamental systems of functions should be no narrower than the interval of the interpolated function.

Representing interpolation polynomials in the form (1) has significant advantages over other forms of representation. For example, there is no need to calculate coefficients $c_k$ interpolation polynomial; The value of the same functions $\varphi_0(t), \varphi_1(t), \ldots, \varphi_N(t)$ can be calculated in advance. Further, representation (1) is linearly dependent on the values of the function $f_j$, $(j = 0,1,\ldots,N)$, which in many cases is essential. Finally, the construction of interpolation polynomials of several variables is greatly simplified, since the fundamental functions of several variables represent the product of the fundamental functions for each variable. The fundamental systems of functions best known till now include the system of fundamental Lagrange interpolation functions, the system of fundamental interpolation trigonometric polynomials, and the Kotelnikov - Shannon system of functions. However, there are other systems of fundamental functions [4]. We restrict ourselves to the fundamental trigonometric polynomials and splines given on uniform grids; It should be noted that trigonometric splines were considered in [5]. We also consider the problem of approximation of functions by the least-squares method and present the systems of fundamental trigonometric polynomials and splines that can be used to construct such approximations.

## The purpose of the work.

Construction of systems of fundamental trigonometric polynomials and splines on uniform grids of two types, used for interpolation of functions and for approximation of functions by the method of least squares.

## The main part.

### System of fundamental trigonometric interpolation polynomials.

Considering the trigonometric functions here and forward, we will, without limiting the generality, assume that $T \equiv 2\pi$. In this segment we will consider uniform grids $\Delta_N^{(I)} = \{t_j^{(I)}\}_{j=1}^N$, ($I = 0, 1$), where

$$t_j^{(0)} = \frac{2\pi}{N}(j-1), \quad t_j^{(1)} = \frac{\pi}{N}(2j-1), \quad N = 2n-1, \quad (n = 1, 2, \ldots).$$

It is known [6] that the trigonometric polynomials fundamental on these grids are fundamental $tm_k(t)$, ($k = 1, 2, \cdots, N$), order $n$ can be written as follows

$$tm_k^{(I)}(t) = \frac{1}{N}\left[1 + 2\sum_{j=1}^n \cos j\left(t - t_k^{(I)}\right)\right].$$

The system of fundamental trigonometric functions generated by grid $\Delta_N^{(I)}$, is the only one in the sense that there is only one system of trigonometric polynomials of order $n$, which satisfies the conditions on this grid

$$tm_k^{(I)}\left(t_j^{(I)}\right) = \begin{cases} 1, & k = j; \\ 0, & k \neq j. \end{cases} \quad (k, j = 1, \ldots, N).$$

Graphs of some fundamental trigonometric polynomials on a grid $\Delta_N^{(0)}$ are shown in Pic.1; note that here and forward we put $N = 9$.

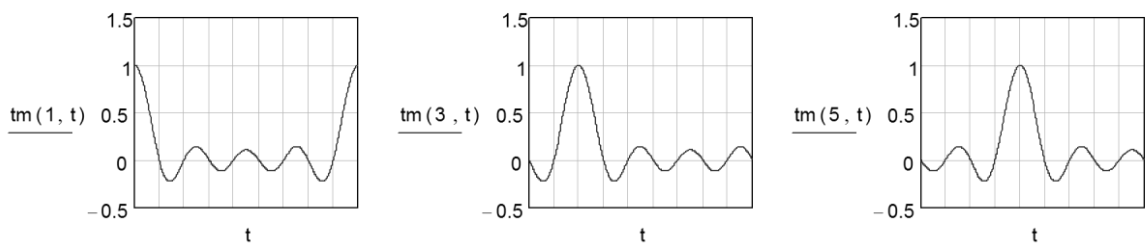

*Pic 1.* Fundamental on the grid $\Delta_N^{(0)}$ trigonometric polynomials $tm_1(t)$, $tm_3(t)$ and $tm_5(t)$.

Graphs of some fundamental trigonometric polynomials on a uniform grid $\Delta_N^{(1)}$ are shown in Pic. 2. Note that the nodes of the grid are $\Delta_N^{(1)}$ are located between the nodes of the grid $\Delta_N^{(0)}$ and are not displayed on the graph.

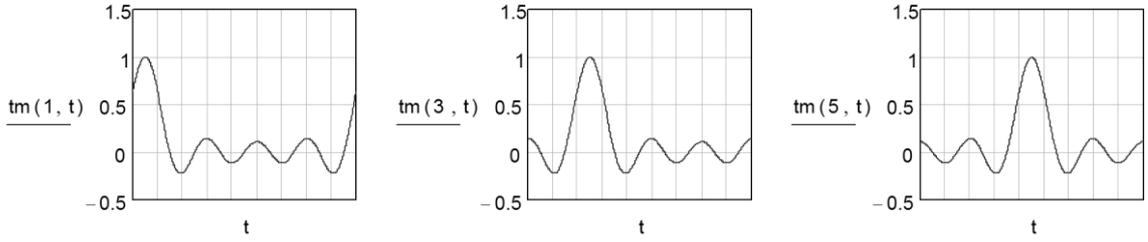

*Pic 2.* Fundamental on the grid $\Delta_N^{(1)}$ trigonometric polynomials $tm_1(t)$, $tm_3(t)$ and $tm_5(t)$.

Using the system of fundamental trigonometric polynomials $tm_k^{(I)}(t)$, $k = 1, 2, \cdots, N$, the interpolation trigonometric polynomial can be written as

$$T_n^{(I)}(t) = \sum_{k=1}^{N} f_k^{(I)} tm_k^{(I)}(t).$$

We note that fundamental trigonometric interpolation polynomials are functions with double orthogonality, that is, they are orthogonal to the norm in understanding the continuous scalar product

$$\sqrt{\frac{N}{2\pi}} \int_0^{2\pi} tm_i^{(I)}(t) tm_j^{(I)}(t) \, dt = \begin{cases} 1, & i = j; \\ 0, & i \neq j, \end{cases} \quad i, j = 1, 2, \ldots, N,$$

and orthogonal in the sense of a discrete scalar product, that is

$$\sum_{k=1}^{N} tm_i^{(I)}(t_k^{(I)}) tm_j^{(I)}(t_k^{(I)}) = \begin{cases} 1, & i = j; \\ 0, & i \neq j. \end{cases} \quad i, j = 1, 2, \ldots, N.$$

**System of fundamental trigonometric interpolation splines.**

Trigonometric interpolation splines on grids $\Delta_N^I$ considered in [5]. In this paper, we will restrict ourselves to only those trigonometric splines that have polynomial analogs.

Fundamental trigonometric interpolation splines can be written in the form [7]

$$ts_j^{(I)}(\sigma, r, N, t) = \frac{1}{N} \left\{ 1 + 2 \sum_{k=1}^{\frac{N-1}{2}} \frac{C_k^{(I)}(\sigma, r, N, j, t)}{H_k^{(I)}(r, N)} \right\};$$

where

$$C_k^{(I)}(\sigma, r, N, j, t) = \sigma_k(r) \cos k(t - t_j^{(I)}) +$$

$$+ \sum_{m=1}^{\infty} (-1)^{mI} \left[ \sigma_{mN+k}(r) \cos(mN + k)(t - t_j^{(I)}) + \sigma_{mN-k}(r) \cos(mN - k)(t - t_j^{(I)}) \right];$$

$$H_k^{(I)}(r, N) = \sigma_k(r) + \sum_{m=1}^{\infty} (-1)^{mI} \left[ \sigma_{mN+k}(r) + \sigma_{mN-k}(r) \right];$$

$$\sigma_k(r, N) = \left[ \frac{\sin(.5hk)}{k} \right]^{1+r}, \quad h = \frac{2\pi}{N}.$$

Graphs of fundamental trigonometric splines $ts_k^{(0)}(r, t)$ at different parameter values $r$ are shown in pic. 3. Note that such trigonometric splines with odd values of this parameter have polynomial analogues - simple polynomial splines of odd degree; for even values of the parameter, the polynomial analogs of these splines are unknown.

Note that forward we will omit the fixed parameters in the graphs; so, for example, since the number of grid nodes $N = 9$ already defined, the graphs do not indicate the dependence from these parameters.

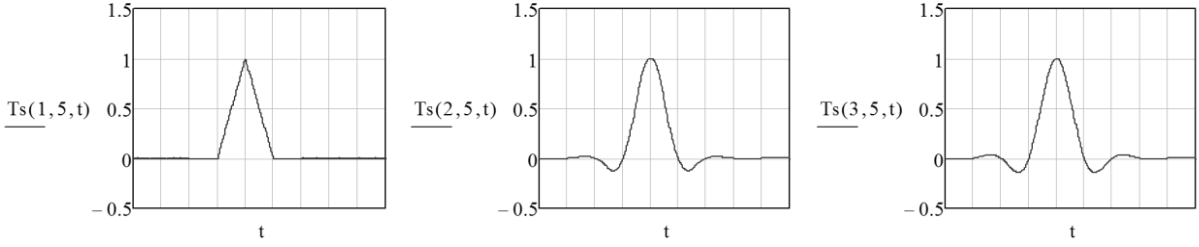

*Pic 3.* Fundamental on the grid $\Delta_N^{(0)}$ trigonometric splines $ts_j^{(0)}(r,t)$ while $j=5$ and parameter values $r=1,2,3$.

Fundamental trigonometric splines $ts_k^{(1)}(\sigma,r,t)$ with odd parameter values $r$ have no polynomial analogues, and for even parameter values, polynomial analogs are simple polynomial splines of even degree. The graphs of such splines at different parameter values are shown in Pic. 4.

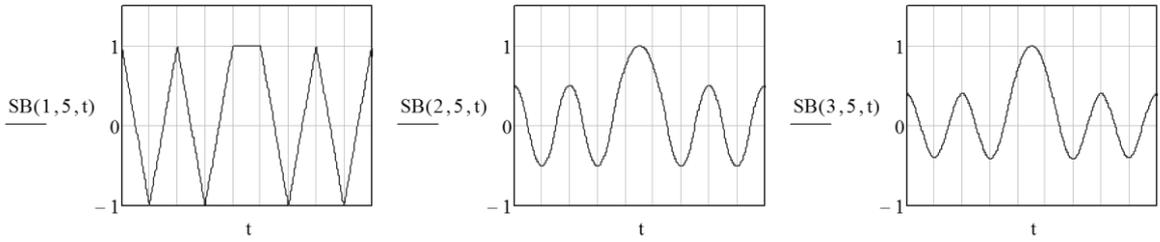

*Pic. 4.* Fundamental on the grid $\Delta_N^{(1)}$ trigonometric splines $ts_j^{(1)}(r,t)$ while $j=5$ and parameter values $r=1,2,3$.

Using fundamental trigonometric interpolation splines $ts_j^{(I)}(r,t)$, interpolation trigonometric spline $St^{(I)}(r,t)$, that interpolates the function $f(t)$ in the nodes of grid $\Delta_N^{(I)}$, can be submitted as

$$St^{(I)}(r,t) = \sum_{k=1}^{N} f_k^{(I)} ts_k^{(I)}(r,t).$$

Fundamental trigonometric splines are orthogonal only in the understanding of a discrete scalar product, that is

$$\sum_{k=1}^{N} ts_i^{(I)}(t_k^{(I)}) ts_j^{(I)}(t_k^{(I)}) = \begin{cases} 1, i = j; \\ 0, i \neq j. \end{cases} \qquad i,j = 1,2,\ldots,N.$$

**Least squares approximation of functions.**

In many cases, instead of interpolating the functions given by their values at the nodes of a grid, it is advisable to use their least squares (LS) approximation. It is known that such a problem is formulated as follows.

As before, on segment $[0,T]$ consider grids $\Delta_N^{(I)}$ and some function $f(t)$, moreover, values are known $f(t_j) = f_j$, $(j = 0,1,\ldots,N)$, of this function in the grid nodes $\Delta_N$. We need to construct a generalized polynomial

$$\Phi_N(t) = c_0\varphi_0(t) + c_1\varphi_1(t) + \ldots + c_N\varphi_N(t)$$

on a system of linearly independent function $\varphi_0(t),\varphi_1(t),\ldots\varphi_N(t)$, which depends from $q$, $(q \leq N)$ parameters $c_0,c_1,\ldots,c_q$, at which the minimum value is reached

$$E_N^{(I)}(f,T_q) = \sum_{j=0}^{N}\left[ f(t_j^{(I)}) - T_q(t_j^{(I)}) \right]^2. \qquad (3)$$

It is known that, together with the classical theory of Fourier series (see, [8]), which considers the continuous case of criterion (3), there is a theory that considers functions on a discrete set of equidistant points; it is clear that in this case appear finite Furge series [2]. In this case, the best approximation for the LS is provided by the segment of finite Fourier series; this fact we will use in constructing the required trigonometric polynomials.

First of all, we note that the interpolation trigonometric polynomial $T_n^{(I)}(t)$, which interpolates the function $f(t)$ in the nodes of grid $\Delta_N^{(I)}$, can be written as

$$T_n^{(I)}(t) = \frac{a_0^{(I)}}{2} + \sum_{k=1}^{n}\left(a_k^{(I)} \cos kt + b_k^{(I)} \sin kt\right), \qquad (4)$$

Where coefficients $a_0^{(I)}, a_k^{(I)}, b_k^{(I)}$, ($k=1,2,...,n$) are calculated by formulas

$$a_0^{(I)} = \tfrac{2}{N}\sum_{j=1}^{N} f_j^{(I)};$$

$$a_k^{(I)} = \tfrac{2}{N}\sum_{j=1}^{N} f_j^{(I)} \cos kt_j^{(I)}; \qquad a_k^{(I)} = \tfrac{2}{N}\sum_{j=1}^{N} f_j^{(I)} \sin kt_j^{(I)}. \qquad (5)$$

Because the system functions $1$, $\cos kt_j^{(I)}$, $\sin kt_j^{(I)}$ is an orthogonal system in the sense of a discrete scalar product [2], then expression (4) can be considered as a finite Fourier series with coefficients (5). It follows from the general theory of Fourier series that partial sums of finite series (4) provide the best approximation for LS in the sense that for anyone $q$, ($q \le n$), minimum values are $E_N(f, T_q)$ is achieved by polynomials

$$T_q^{(I)}(t) = \frac{a_0^{(I)}}{2} + \sum_{k=1}^{q}\left(a_k^{(I)} \cos kt + b_k^{(I)} \sin kt\right) \qquad (6)$$

with coefficients (5). It is clear that when $q = n$ value is $E_N(f, T_N) = 0$.

In what follows, to reduce the fundamental trigonometric polynomials and the fundamental trigonometric splines approximating the function on a discrete set of points in the LS, we will refer to the fundamental trigonometric LS polynomials and the fundamental trigonometric LS splines, respectively.

The problem of constructing systems of fundamental trigonometric LS polynomials and fundamental trigonometric LS splines is still relevant. Let's take a closer look at this problem.

**System of fundamental trigonometric LS polynomials.**

We will build such a system as follows: substituting expressions (5) for the coefficients of the polynomial $T_q^{(I)}(t)$ in (6), after changing the order of summation we get

$$T_q^{(I)}(t) = \frac{1}{2}\frac{2}{N}\sum_{j=1}^{N} f_j^{(I)} +$$

$$+ \sum_{k=1}^{q}\left[\frac{2}{N}\sum_{j=1}^{N} f_j^{(I)} \cos kt_j^{(I)} \cos kt + \frac{2}{N}\sum_{j=1}^{N} f_j^{(I)} \sin kt_j^{(I)} \sin kt\right] =$$

$$= \sum_{j=1}^{N} f_j^{(I)}\left\{\frac{1}{N} + \frac{2}{N}\sum_{k=1}^{q}\left[\cos kt_j^{(I)} \cos kt + \sin kt_j^{(I)} \sin kt\right]\right\} =$$

$$= \sum_{j=1}^{M} f_j^{(I)}\left\{\frac{1}{N}\left[1 + 2\sum_{k=1}^{q}\cos k\left(t - t_j^{(I)}\right)\right]\right\}.$$

This expression can be represented as

$$T_n^{(I)}(t) = \sum_{j=1}^{N} f_j^{(I)}\, \varphi_{j,q}^{(I)}(t), \qquad (7)$$

where

$$\varphi_{j,q}^{(I)}(t) = \frac{1}{N}\left[1 + 2\sum_{k=1}^{q}\cos k\left(t - t_j^{(I)}\right)\right], \qquad (8)$$

$$j = 1,\ldots,N\,;\ q \le n.$$

Therefore, the system of fundamental trigonometric LS polynomials consists of polynomials $\varphi_{j,q}^{(I)}(t)$, defined by expression (8). Graphs of some LS polynomials $\varphi_{j,q}^{(I)}(t)$ on grids $\Delta_N^{(I)}$ when fixed $j$ and $N$ for different parameter values $q$ is given on Pic. 5,6.

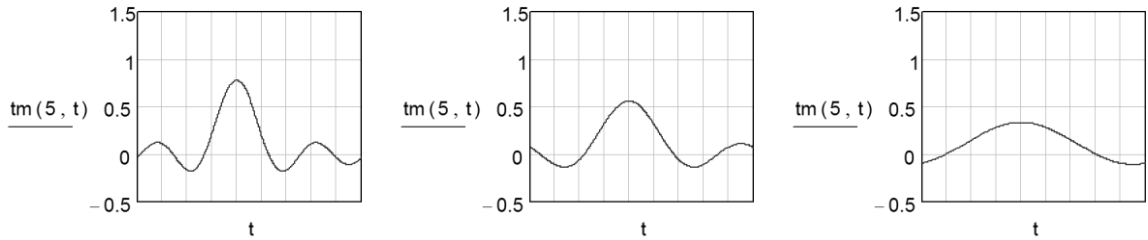

*Pic.5.* Fundamental trigonometric LS polynomials
on the grid $\Delta_N^{(0)}$; $j=5$, $q=3,2,1$.

**System of fundamental trigonometric LS splines.**

A certain disadvantage of the polynomial order $m$ the best standard deviation approximation, based on the values of a function on a discrete set of equidistant points, is its analyticity, which makes this polynomial not very convenient for approximating functions with a small smoothness. To approximate such functions, it is more convenient to use trigonometric splines whose values at the nodes of the grid coincide with the values of the trigonometric LS polynomial of the order $q$, constructed by the values of a function on a discrete set of equidistant points.

It is easy to show that trigonometric spline is of order $q$, the values of which at the nodes of the grid coincide with the values of a polynomial of discrete approximation of the same order, built on the values of the function on a discrete set of equidistant points, can be represented as $Ts_{r,q}^{(I)}(x) = \sum_{j=1}^{N} f_j ts_{j,q}^{(I)}(r,x)$,

where

$$ts_{j,q}^{(I)}(r,t) = \frac{1}{N}\left\{1 + 2\sum_{k=1}^{q} \frac{C_k^{(I)}(r,N,j,t)}{H_k^{(I)}(r,N)}\right\},$$

And functions $C_k^{(I)}(\sigma,r,N,j,t)$ and constans $H_k^{(I)}(r,N)$ and $\sigma_k(r,N)$ are defined in the same way as before. Trigonometric splines $ts_{j,q}^{(I)}(r,x)$ we will call fundamental trigonometric LS splines. Graphs of some fundamental trigonometric LS splines $ts_{j,q}^{(I)}(r,x)$ on grids $\Delta_N^{(I)}$ when fixed $j$ for different parameter values $q$ and $r$ are given on Pic. 6-9.

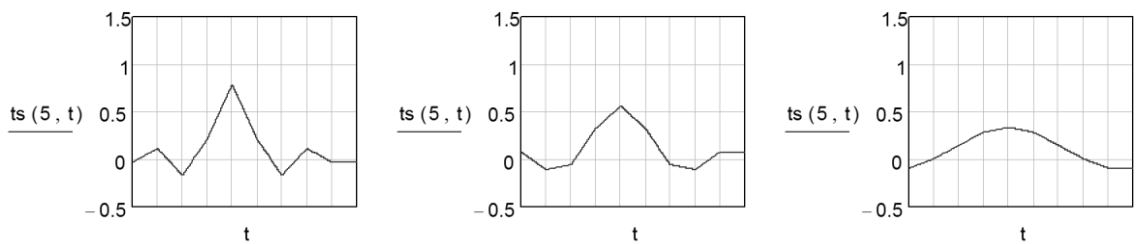

*Pic.6.* Fundamental trigonometric LS splines
on the grid $\Delta_N^{(0)}$; $r=1$; $q=3,2,1$.

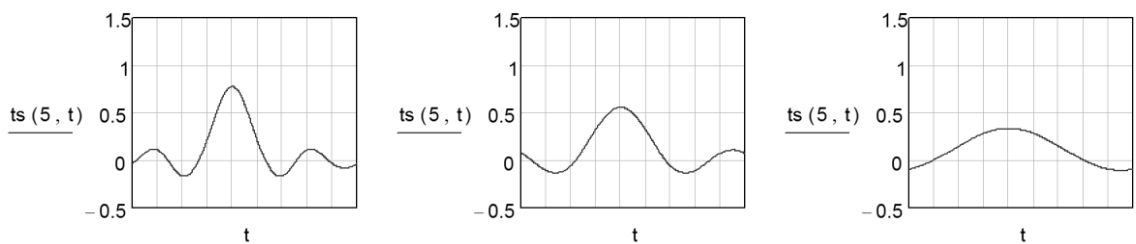

*Pic.7.* Fundamental trigonometric LS splines
on the grid $\Delta_N^{(0)}$; $r=3$; $q=3,2,1$

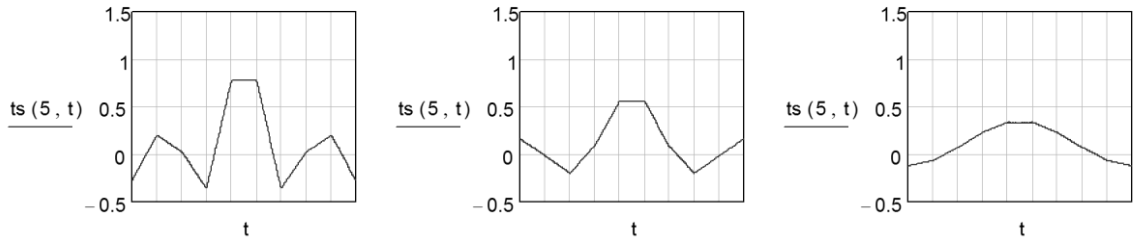

*Pic.8.* Fundamental trigonometric LS splines
on the grid $\Delta_N^{(1)}$; $r=1$; $q=3,2,1$.

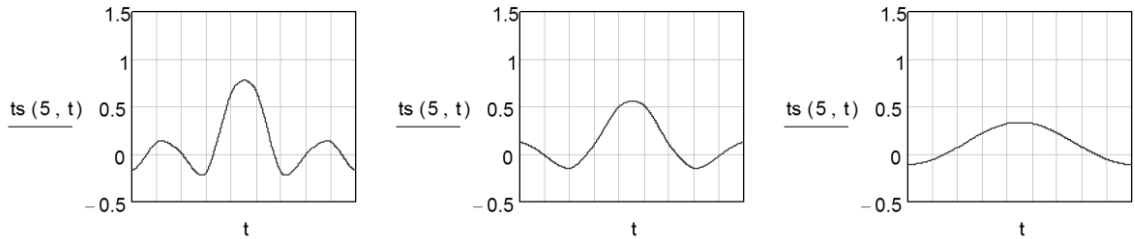

*Pic.9.* Fundamental trigonometric LS splines
on the grid $\Delta_N^{(1)}$; $r=2$; $q=3,2,1$.

Note that as shown in Pic. 6,8, trigonometric LS spline on the grids $\Delta_N^{(I)}$ at parameter value $r=1$, consists of polygons that are stitched at the nodes of the grid $\Delta_N^{(0)}$.

It should be noted that fundamental trigonometric LS polynomials and LS splines lose their orthogonality properties in a continuous and discrete sense.

## Conclusions

1. Fundamental trigonometric interpolation polynomials are constructed on uniform grids $\Delta_N^{(I)}$, ($I=0,1$).
2. Fundamental trigonometric interpolation splines are constructed on uniform grids $\Delta_N^{(I)}$, ($I=0,1$).
3. Fundamental trigonometric LS polynomials are constructed on uniform grids $\Delta_N^{(I)}$, ($I=0,1$), which are used to approximate functions using the least squares method.
4. Fundamental trigonometric LS splines are constructed on uniform grids $\Delta_N^{(I)}$, ($I=0,1$), which are used to approximate functions using the least squares method.
5. Of course, the construction of fundamental trigonometric interpolation polynomials and splines on uniform grids $\Delta_N^{(I)}$, ($I=0,1$) is an interesting task, for other cases considered in [7]; similarly, it is of interest to construct fundamental trigonometric LS polynomials and LS splines on uniform grids $\Delta_N^{(I)}$, ($I=0,1$), for the same cases.
6. The following fundamental trigonometric polynomials and splines can be recommended for use in many problems related to approximation of functions, in particular mathematical modeling of signals and processes.

## List of references

1. Berezin I. S., Zhidkov N. P. Computing Methods. Vol. I und II. Oxford ,1965. Pergamon Press.
2. Hamming R.W. Numerical Methods for Scientists and Engineers. –MC Grow-Hill Book Company, INC. - 1962.
3. R. S. Guter, L. D. Kudryavtsev, B. M. Levitan, Elements of the Theory of Functions - Translation edited by I. N. Sneddon.- Pergamon, 1966.

4. Denysiuk V.P. Fundamental Functions and Trigonometrical Splines – Kiev.: ПАТ "Випол", 2015. (Ukraine)
5. Denysiuk V. P. Generalized Trigonometric Functions and Their Applications - IOSR Journal of Mathematics (IOSR-JM), Volume 14, Issue 6 Ver. I (Nov - Dec 2018), PP 19-25.
6. Natanson I.P. Constructive Theory of Functions – Gosudarsnvennoe Izdatel'stvo Techniko - Teoreticheskoi Literatury, Moskva, Leningrad, 1949.
7. Kolmogorov A.N., Fomin S.V. Elements of the Theory of Functions and Functional Analysis – Vol. 1 and 2. Dover Publications, INC. Mineola, New York.
8. Denysiuk V.P. About the classification of trigonometric splines - ArXiv:1910.00830